\documentclass[11pt,a4paper]{article}
\usepackage[latin1]{inputenc} 
\usepackage[T1]{fontenc}
\usepackage{amsmath, amssymb, dsfont,eufrak}
\usepackage{graphics, pstricks, epsfig, psfig}
\usepackage[frenchb]{babel}
\usepackage[all]{xy}
\usepackage{float}
\usepackage{multicol}
\usepackage{here}
\usepackage{array}

\def\color[#1]#2{}

\newtheorem{theo}{Th{\'e}or{\`e}me}

\newtheorem{lemme}{Lemme}
         \newenvironment{demo}{\textbf{D{\'e}monstration} :}{}

\newtheorem{propo}{Proposition}

\title{
Existence d'une courbe sur $\mathbb{F}_3$ de genre $5$ avec
  $13$ points rationnels}
\author{Christophe Ritzenthaler}

\begin{document}
\maketitle

\begin{abstract}
Let $N_q(g)$ the maximal number of points on a genus $g$ curve over
$\mathbb{F}_q$. We prove that $N_3(5)=13$.
\end{abstract}

\section{Introduction}
{\def\thefootnote{}\footnotetext{\hskip-1.8em Univ. Paris VII, laboratoire de
    th{\'e}orie des nombres. \\  E-mail : ritzenth@math.jussieu.fr}}%
Soit $N_q(g)$ le nombre maximal de points rationnels sur $\mathbb{F}_q$  
pour une courbe lisse de genre $g$ sur $\mathbb{F}_q$.
Gr{\^a}ce aux bornes de Weil on sait que $N_q(g)$ est 
 inf{\'e}rieur ou {\'e}gal {\`a} $q+1+2g \sqrt{q}$. Ces
bornes (m{\^e}me raffin{\'e}es par la m{\'e}thode d'Oesterl{\'e}) ne sont pas 
optimales et la d{\'e}termination de $N_q(g)$ pour $g>2$ reste incompl{\`e}te. \\
Le cas qui nous pr{\'e}occupe ici est celui des courbes de genre $5$ sur
$\mathbb{F}_3$. Les majorations explicites d'Oesterl{\'e} donne un
nombre de points inf{\'e}rieur ou {\'e}gal {\`a} $14$, et 
Kristin Lauter a montr{\'e} l'in{\'e}galit{\'e} stricte (cf. \cite{kristin}). 
D'autre part on
connaissait l'existence de courbes avec $12$ points rationnels
construites par des rev{\^e}tements successifs (cf. \cite{niederreiter}). Dans
la pr{\'e}sente note nous nous proposons donc de combler la lacune
existante en donnant explicitement
une courbe de genre $5$ avec $13$ points rationnels comme
l'intersection de $3$ quadriques dans $\mathbb{P}^4$. Cette courbe est
de plus exceptionnelle pour une autre raison : elle constitue, {\`a} ma
connaissance, le premier exemple d'une courbe avec un nombre de points
maximum qui est rev{\^e}tement 
non galoisien d'une courbe elliptique. Nous donnons de plus
explicitement ce dernier.

\section{R{\'e}sultats}
\begin{propo}
La courbe $C$ d{\'e}finie par 
$$\left\{
\begin{array}{l}
q_1 =2 x_1 x_2+x_3 x_2+x_3^2-x_4^2 \\
q_2 =x_5 x_1-x_4 x_2 \\
q_3 =x_1^2+x_1 x_2-x_3^2+x_5^2 \\
\end{array}
\right.$$
est une courbe de genre $5$ qui poss{\`e}de $13$ points sur
$\mathbb{F}_3$. 
\end{propo}
\begin{propo}
\begin{enumerate}
\item $\textrm{Aut}_{\overline{\mathbb{F}}_3}(C)=<\omega> \simeq
\mathbb{Z}/2\mathbb{Z}$.  
\item $D=C/<\omega>$ a pour {\'e}quation
$$x_2^2(x_3 x_2+x_3^2-x_1 x_2)+(x_1^2+ x_1 x_2-x_3^2)
x_1^2=0.$$
\item $C$ est rev{\^e}tement non galoisien de la courbe elliptique $E :
  y^2=x^3-x+1$.  De plus si on prend comme mod{\`e}le plan pour $C$
$$x^4+x^3 y^3-x^2-x y^5+y^5+2 y=0$$
le rev{\^e}tement est donn{\'e} par $(x:y:1) \mapsto (x':y':1)$ avec 
$$\begin{cases}
x'=\frac{-x^3-x^2y^3-x^2y^2+x^2+xy^5+xy^4-xy^3-xy^2+xy-y^6+y^5+y^4+y^3-y^2}{(y+1)
  (y-1)^2 (y^3-y^2+y+1)} \\ 
y'= \frac{-x^3 y-x^2 y^5+x^2y^4-x^2y+xy^5+xy^4+xy^3-xy^2-xy+y^8+y^7+y^4-y^3-y^2-y+1
}{(y+1)^2 (y-1)^3 (y^3-y^2+y+1)}
\end{cases}
$$
\end{enumerate}
\end{propo}

\section{D{\'e}monstration}
La courbe a {\'e}t{\'e} construite par une recherche exhaustive des
sextiques planes de genre $5$ passant par les $13$ points rationnels
du plan projectif
(soit $3^{15}$ possibilit{\'e}s). Le plongement canonique permet alors
d'obtenir un mod{\`e}le lisse comme intersection de $3$ quadriques dans
$\mathbb{P}^4$.\\
 
D'autre part le polyn{\^o}me caract{\'e}ristique de $C$ sur $\mathbb{F}_3$
 se factorise en 
\begin{equation} \label{polycarac}
(T^2+2T+3)(T^2+3T+3)(T^2+3)(T^4+4T^3+8T^2+12T+9).
\end{equation} 
De plus sur $\mathbb{F}_{3^{24}}$, le polyn{\^o}me caract{\'e}ristique se
scinde en
$$(T-3^{12})^4 (T^2+629918 T +3^{24})^3.$$
Sur $\mathbb{F}_{3^{24}}$ la jacobienne de la courbe $C$ est donc
isog{\`e}ne {\`a} $E^2 \times F^3$, avec $E$ et $F$ des courbes
elliptiques qui sont absolument non-isog{\`e}nes (par exemple parce que
la premi{\`e}re est supersinguli{\`e}re et pas l'autre). En particulier sur
$\mathbb{F}_3$, $C$ est rev{\^e}tement de trois courbes elliptiques.\\

Pour montrer que $C$ n'est rev{\^e}tement galoisien d'aucune de ces
courbes elliptiques, nous allons d{\'e}terminer 
le groupe des automorphismes de $C$.
On constate que $$\omega : (x_1,x_2,x_3,x_4,x_5) \mapsto
(x_1,x_2,x_3,-x_4,-x_5)$$ est un automorphisme de la courbe. De plus
$C/<\omega>$ est une courbe de genre $2$ (par la formule d'Hurwitz) 
qu'on peut obtenir
par l'{\'e}limination des variables $x_4$ et $x_5$ sous la forme :
$$2 y^3 x+ y^3+y^2+x^4+x^3 y-x^2.$$

Montrons maintenant que cet automorphisme est le seul qui soit non
trivial, on aura alors que $C$ est rev{\^e}tement non galoisien d'une
courbe elliptique. C'est en
fait une cons{\'e}quence d'un th{\'e}or{\`e}me plus g{\'e}n{\'e}ral d{\`u}
{\`a} Beauville \cite[prop.~6.9]{beau} et qui donne exactement le groupe
des automorphismes de $C$ en fonction de ceux de la quintique
d{\'e}finie ci-dessous. Mais nous avons besoin de quelque chose de moins 
pr{\'e}cis et on obtient la d{\'e}monstration {\'e}l{\'e}mentaire ci-dessous. \\
Puisque $C$ est donn{\'e}e sous forme canonique, tous les
automorphismes de la courbe sont lin{\'e}aires. Soit $\psi$ un 
  automorphisme de $\mathbb{P}^4$. C'est un automorphisme de
$C$ si et seulement si la 
matrice $M \in \rm{PGL}_5(\overline{\mathbb{F}_3})$ qui le
repr{\'e}sente  est telle que
$q_i(M v)=0$ pour $i=1,2,3$ et quelque soit $^t
v=(x_1,x_2,x_3,x_4,x_5) \in C$.\\
Mais si $M$ repr{\'e}sente un automorphisme de $C$ alors 
si on note $Q_i$ les matrices des formes quadratiques
$q_i$, $^t M Q_i M$ est une 
quadrique contenant $C$. Elle est donc combinaison
lin{\'e}aire des $Q_i$.\\
On consid{\`e}re alors l'ensemble $S$ des $(x,y,z) \in
\mathbb{P}_2$ tels que $\rm{det}(x \, Q_1+ y \, Q_2 + z \, Q_3)=0$. C'est une
quintique lisse d'{\'e}quation 
$$-x^3+y^2 x-y^2-x^4+x+x^3 y^2-y^4 x+y^4=0$$
qui ne poss{\`e}de qu'un seul automorphisme $$\varphi : (x,y,z)
\mapsto (x,-y,z).$$ On peut d{\'e}finir un morphisme de groupe $\mu$ de
$\rm{Aut}(C) \subset \rm{PGL}_5(\overline{\mathbb{F}_3})$
dans $\rm{Aut}(S)$ : Si $M$ est un 
automorphisme de la courbe $C$ et si 
$$\left\{
\begin{array}{l}
^t M Q_1 M =a_1 Q_1+ b_1 Q_2 + c_1 Q_3 \\
^t M Q_2 M =a_2 Q_1+ b_2 Q_2 + c_2 Q_3 \\
^t M Q_3 M =a_3 Q_1+b_3 Q_2+ c_3 Q_3 \\
\end{array}
\right.$$
on a alors un automorphisme de $S$ donn{\'e}e par $(x,y,z) \mapsto (a_1
x+a_2 y+a_3 z, b_1x+b_2y+b_3z,c_1x+c_2y+c_3z)$. 
De plus ce morphisme envoie $\omega$ sur $\varphi$.

Il suffit de montrer que $\mu$ est injectif. 
La quintique {\'e}tant non
  singuli{\`e}re, la quadrique singuli{\`e}re $x Q_1+y Q_2+z Q_3$ 
associ{\'e}e {\`a} un point de
  la courbe $(x,y,z)$ est de rang $4$ et poss{\`e}de donc un unique point
  singulier. Soit
  $s : S \rightarrow \mathbb{P}^4$ qui associe {\`a} un 
  point de la quintique le point singulier de la quadrique
  correspondante. Soit $M$ un automorphisme de $C$ qui se
  r{\'e}duit sur l'identit{\'e} de $S$. 
  L'action de $M$ en tant 
  qu'automorphisme de $\mathbb{P}^4$ sur $s(S)$ est la
  m{\^e}me que celle induite par $s^*(\mu^*(M))$ qui est dans ce cas
  l'identit{\'e}. Pour montrer qu'on a alors $M=\mathrm{Id}$ il
  suffit donc de montrer que les points de $s(S)$ ne sont pas
  contenus dans un hyperplan. C'est en fait une cons{\'e}quence du lemme
  suivant :
\begin{lemme}
Soit $F$ et $G$ deux matrices de  formes quadratiques non
d{\'e}g{\'e}n{\'e}r{\'e}es d'un espace
vectoriel $E$ sur un corps $k$ de caract{\'e}ristique diff{\'e}rente de
$2$ 
 telles que $P(t)=\mathrm{det}(G-t F)=0$ n'ait que des
racines de multiplicit{\'e} $1$. Alors ces deux formes quadratiques sont
simultan{\'e}ment diagonalisables.
\end{lemme}
Montrons tout d'abord comment ce lemme permet de conclure.
Soit $l$ une droite transverse {\`a} $S$ et $p_0,\ldots,p_4$ les points
d'intersection. On consid{\`e}re le pinceau de quadriques d{\'e}fini par
$l$ : il est engendr{\'e} par deux quadriques non singuli{\`e}res
d'{\'e}quations $F$ et $G$ telles que $\mathrm{det}(F-tG)=0$ n'a que des
  racines de multiplicit{\'e}s $1$ (puisque $l$ est transverse). On peut
  alors appliquer le lemme : dans 
une base qui diagonalise simultan{\'e}ment les deux
  quadriques on {\'e}crit $F=\sum X_i^2$ et $G=\sum \alpha_i X_i^2$
  $\alpha_i \neq \alpha_j$. Avec ces coordonn{\'e}es les images $s(p_i)$
  sont alors tout simplement les points
  $(1,0,0,0,0),(0,1,0,0,0),\ldots,(0,0,0,0,1)$ qui ne sont
  {\'e}videmment pas dans un m{\^e}me hyperplan, d'o{\`u} le r{\'e}sultat.\\

\begin{demo}
Soit $n$ la dimension de $E$.
Soient $\lambda_i,v_i$ ($\lambda_i \neq 0$ et $v_i \neq 0$) les $n$
scalaires et vecteurs tels que $F v_i=\lambda_i G_i v_i$. Nous allons
montrer que les $v_i$ sont une base orthogonale pour $F$ et $G$.
On a
\begin{eqnarray}
^t v_j F v_i  &= &  \lambda_i \;  ^t v_j  G v_i   \\
 &  =& ^t v_i F v_j  \; \textrm{par sym{\'e}trie de } \; F \\ \label{lig1}
 & = & \lambda_j \;  ^t v_i  G v_j  \\
 & = & \lambda_j \; ^t v_j  G v_i \;  \textrm{par sym{\'e}trie de} \; G \label{lig4}
\end{eqnarray}
Par hypoth{\`e}se, les $\lambda_i$ sont tous distincts on a donc par
{\'e}galit{\'e} de \ref{lig1} et \ref{lig4}  que $^t v_j F v_i={^t
v_j} G v_i =0$.\\ 
\end{demo}

Nous allons donner explicitement un rev{\^e}tement de $C$ sur une courbe
elliptique.
On a le th{\'e}or{\`e}me suivant
\begin{theo} {\rm \cite{lauter2}}
Soit $C$ une courbe sur $\mathbb{F}_q$ dont la jacobienne est
isog{\`e}ne {\`a} un produit $\Lambda \times E$ avec $E$ une courbe
elliptique. Soit $r$ le r{\'e}sultant des polyn{\^o}mes minimaux de la
restriction de $F+V$ {\`a} $E$ et $\Lambda$ o{\`u} $F$ est l'endomorphisme
de Frobenius de $\textrm{Jac}(C)$ et $V$ son dual. Alors il existe une
courbe elliptique $E'$ isog{\`e}ne {\`a} $E$ et une application de $C$
dans $E'$ dont le degr{\'e} divise $r$.
\end{theo}

\begin{figure} 
\caption{Ramification}
\label{ramimi}
      \leavevmode
      \begin{center}
      \input{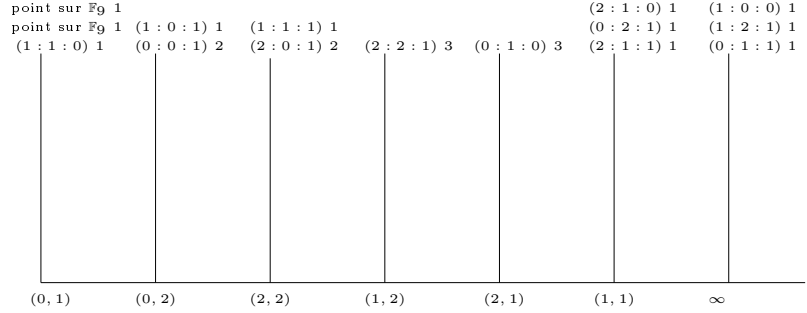}
      \end{center}
\end{figure}

On applique ce th{\'e}or{\`e}me au facteur $T^2+3 T+3$ de (\ref{polycarac}).
 Il existe donc un rev{\^e}tement de degr{\'e} $3$ de $C$
vers une courbe $E$ avec $7$ points sur $\mathbb{F}_3$.
A isomorphisme pr{\`e}s cette courbe est unique d'{\'e}quation $E :
y^2=x^3-x+1$. \\
Pour expliciter le rev{\^e}tement, nous proc{\'e}dons comme suit 
\begin{itemize}
\item Au dessus d'au moins un point rationnel de $E$ il
  existe trois points rationnels de $C$ ({\'e}ventuellement non
  distincts). 
Quitte {\`a} effectuer une
  translation on peut supposer que ce point est l'origine de la courbe
  elliptique.
\item Pour chacune des ${13 \choose 1}+{13 \choose 2}+{13 \choose
    3}=377$ 
possibilit{\'e}s, on
  consid{\`e}re alors le diviseur $D$ de degr{\'e} $3$ au dessus de l'origine. Si $f : C
  \rightarrow E$ est le rev{\^e}tement 
 et si $g : E  \rightarrow \mathbb{P}^1$ est
  l'application $(x:y:z) \mapsto (x:z)$ alors $g \circ f \in
  \mathcal{L}(2 D)$. Le th{\'e}or{\`e}me de Clifford montre que $l(2 D) \leq
  2$ et par Riemann-Roch on a que $l(D)=2$. 
\item Soit $\phi \in
  \mathcal{L}(2 D)$ non constante. Pour $\psi_x=\phi,\phi+1,\phi-1$ on
  calcule $\psi^3-\psi+1$. Si cette fonction est le carr{\'e} d'une
  autre fonction $\psi_y$ alors on a le rev{\^e}tement donn{\'e} par
  $p \mapsto (\psi_x(p):\psi_y(p):1)$.
\end{itemize} 
Gr{\^a}ce {\`a} MAGMA, on r{\'e}alise rapidement ces calculs. Deux mod{\`e}les
plans paraissent particuli{\`e}rement int{\'e}ressants pour $C$ : le
premier $$x^4y^2+
x^2 y^4+y^6+x^2 - y^2 +x^3 +x y^2 +1+x -x^4=0$$ est un mod{\`e}le pour lequel
l'involution est $y \mapsto -y$. Le second $$
x^4+x^3 y^3-x^2-x y^5+y^5+2 y=0$$ poss{\`e}de $13$ points rationnels (tous
les points de $\mathbb{P}^2(\mathbb{F}_3)$). Pour ce second mod{\`e}le,
un rev{\^e}tement est donn{\'e} par $(x:y:1) \mapsto (x':y':1)$ avec 
$$\begin{cases}
x'=\frac{-x^3-x^2y^3-x^2y^2+x^2+xy^5+xy^4-xy^3-xy^2+xy-y^6+y^5+y^4+y^3-y^2}{(y+1)
  (y-1)^2 (y^3-y^2+y+1)} \\ 
y'= \frac{-x^3 y-x^2 y^5+x^2y^4-x^2y+xy^5+xy^4+xy^3-xy^2-xy+y^8+y^7+y^4-y^3-y^2-y+1
}{(y+1)^2 (y-1)^3 (y^3-y^2+y+1)}
\end{cases}
$$
La figure \ref{ramimi} montre la ramification de $f$
(le nombre derri{\`e}re le point est l'indice de ramification. 
On constate en particulier que le rev{\^e}tement est
sauvagement ramifi{\'e}).


\begin{thebibliography}{}
\bibitem{beau} A. Beauville : Vari{\'e}t{\'e}s de Prym et Jacobiennes
  interm{\'e}diaires. Ann. Scient. {\'E}c. Norm. Sup. $4^e$ s{\'e}rie,
  t. 10 (1977), 309-391
\bibitem{lauter2} E.W. Howe \& K. Lauter : Improved upper bounds for
  the number of points on curves over finite fields,
  ArXiv:math.NT/0207101 v5, (2002). 
\bibitem{kristin} K. Lauter : Non-existence of a curve over
    $\mathbb{F}_3$ of genus $5$ with $14$ rational points. Proc. AMS
    {\bf 128} (1999), 369-374
\bibitem{niederreiter} H. Niederreiter, C.P. Xing : Cyclotomic
  function fields, Hilbert class fields and global function fields
  with many rational places, Acta. Arithm. {\bf 79} (1997), 59-76.
\end{thebibliography}
\end{document}